\documentclass[a4paper,11pt,reqno]{amsart}

\usepackage{amssymb}

\newtheorem{theorem}[equation]{Theorem}
\newtheorem{lemma}[equation]{Lemma}
\newtheorem{corollary}[equation]{Corollary}
\newtheorem{proposition}[equation]{Proposition}

\numberwithin{equation}{section}

\theoremstyle{definition}

\newtheorem*{example*}{Example}

\newtheorem{remark}[equation]{Remark}
\newtheorem*{remark*}{Remark}

\newcommand{\bN}{{\mathbb N}}
\newcommand{\bZ}{{\mathbb Z}}

\newcommand{\frg}{{\mathfrak g}}

\newcommand{\frgtetra}{{{\mathfrak g}_{\boxtimes}}}
\newcommand{\frgtetrahat}{{{\mathfrak g}_{\widehat{\boxtimes}}}}

\newcommand{\frf}{{\mathfrak f}}

\newcommand{\frt}{{\mathfrak t}}
\newcommand{\frs}{{\mathfrak s}}
\newcommand{\fra}{{\mathfrak a}}

\newcommand{\frl}{{\mathfrak l}}

\newcommand{\calS}{{\mathcal S}}

\newcommand{\calA}{{\mathcal A}}
\newcommand{\calM}{{\mathcal M}}

\newcommand{\calI}{{\mathcal I}}
\newcommand{\calO}{{\mathcal O}}

\providecommand{\espan}[1]{\text{span}\left\{ #1\right\}}

   \DeclareMathOperator{\lrt}{\mathfrak{lrt}}

  \DeclareMathOperator{\frsldos}{{{\mathfrak{sl}_2}}}

 \DeclareMathOperator{\frgl}{{\mathfrak{gl}}}

 \DeclareMathOperator{\ad}{ad}

 \DeclareMathOperator{\End}{End}

 \DeclareMathOperator{\Aut}{Aut}

\newenvironment{romanenumerate}
 {\begin{enumerate}
 
 }{\end{enumerate}}

\newenvironment{romanprimeenumerate}
 {\begin{enumerate}
 }{\end{enumerate}}

\newenvironment{alphaenumerate}
 {\begin{enumerate}
 
 }{\end{enumerate}}

\begin{document}

\title{The $S_4$-action
on the Tetrahedron algebra}

\author[Alberto Elduque]{Alberto Elduque$^{\star}$}
 \thanks{$^{\star}$ Supported by the Spanish Ministerio de
 Educaci\'on y Ciencia
 and FEDER (MTM 2004-081159-C04-02) and by the
Diputaci\'on General de Arag\'on (Grupo de Investigaci\'on de
\'Algebra)}
 \address{Departamento de Matem\'aticas, Universidad de
Zaragoza, 50009 Zaragoza, Spain}
 \email{elduque@unizar.es}

\date{September 1, 2006}



\begin{abstract}
The action of the symmetric group $S_4$ on the Tetrahedron algebra,
introduced by Hartwig and Terwilliger \cite{Paul}, is studied. This
action gives a grading of the algebra which is related to its
decomposition in \cite{Paul} into a direct sum of three subalgebras
isomorphic to the Onsager algebra. The ideals of both the
Tetrahedron algebra and the Onsager algebra are determined.
\end{abstract}

\maketitle


\section*{Introduction}

The Tetrahedron algebra $\frgtetra$ has been defined in
\cite[Definition 1.1]{Paul}, in connection with the so called
Onsager algebra introduced in \cite{Onsager}, in which the free
energy of the two dimensional Ising model was computed. Since then
it has been investigated by physicists and mathematicians in
connection with solvable lattice models, representation theory,
Kac-Moody Lie algebras, tridiagonal pairs and partially orthogonal
polynomials (see \cite{Paul} and the references there in).

One of the main results (Theorem 11.5) in \cite{Paul} shows that
$\frgtetra$ is isomorphic to the three point $\frsldos$ loop algebra
$\frsldos\otimes \calA$, where $\calA$ is the algebra
$k[t,t^{-1},(1-t)^{-1}]$ (a subalgebra of the field of rational
functions on the indeterminate $t$).

The Tetrahedron algebra is endowed with an action of the symmetric
group $S_4$ by automorphisms, and the result above is used in
\cite[Theorem 11.6]{Paul} to show that $\frgtetra$ is the direct sum
$\frgtetra=\Omega\oplus\Omega'\oplus\Omega''$ of a subalgebra
$\Omega$, which is isomorphic to the Onsager algebra, and its images
under the action of the cycle $(123)$. Hence, $\frgtetra$ is a
direct sum of three subalgebras (not ideals!) which are isomorphic
to the Onsager algebra.

On the other hand, like any other Lie algebra endowed with an action
of $S_4$ by automorphisms, $\frgtetra$ is $\bZ_2\times\bZ_2$-graded
by the action of Klein's $4$ group (see \cite{EO}):
\begin{equation}\label{eq:Z2Z2}
\frgtetra=(\frgtetra)_{(\bar 0,\bar 0)} \oplus
 (\frgtetra)_{(\bar 1,\bar 0)} \oplus
 (\frgtetra)_{(\bar 0,\bar 1)} \oplus
 (\frgtetra)_{(\bar 1,\bar 1)}.
\end{equation}
In \cite[Problem 13.4]{Paul}, the authors pose the problem of
showing that $(\frgtetra)_{(\bar 0,\bar 0)}=0$, of finding a basis
for each of the subspaces $(\frgtetra)_{(\bar 1,\bar 0)}$,
$(\frgtetra)_{(\bar 0,\bar 1)}$ and $(\frgtetra)_{(\bar 1,\bar 1)}$
(which become abelian subalgebras), and of relating this
$\bZ_2\times\bZ_2$-grading to the previous decomposition
$\frgtetra=\Omega\oplus\Omega'\oplus\Omega''$.

\smallskip

These notes give a solution to this problem. A key to this solution
will be the use of a very suitable basis of $\frsldos\otimes \calA$,
as a Lie algebra over $\calA$, which is also useful in simplifying
some of the arguments in \cite{Paul}.

\smallskip

All the algebras considered will be defined over a ground field $k$
of characteristic $\ne 2$. Unadorned tensor products $\otimes$ will
be considered over $k$.

In the next section, the isomorphism between $\frgtetra$ and
$\frsldos\otimes\calA$, as algebras with an action of the symmetric
group $S_4$, given in \cite{Paul}, will be reviewed and proved in a
simplified way, by using the suitable $\calA$-basis mentioned above.
The action of $S_4$ on $\frgtetra$ translates into an action of
Klein's $4$ group as $\calA$-automorphisms of $\frsldos\otimes\calA$
plus an action of $S_3$ on both $\frsldos$ and $\calA$. Section 2
will be devoted to solve \cite[Problem 13.4]{Paul}, and the normal
Lie related triple algebra which is associated to the action of
$S_4$ on $\frgtetra$, as shown in \cite{EO}, will be found in
Section 3. This will highlight a general construction of normal Lie
related triple algebras defined on certain commutative associative
algebras endowed with an action of the symmetric group $S_3$.
Section 4 will be devoted to solve \cite[Problem 13.3]{Paul}, which
asks for the ideals of the Tetrahedron algebra, in terms of this
$\calA$-basis used throughout. A general result on ideals of some
Lie algebras will be given, and the ideals of the Onsager algebra
will be determined too. Section 5 will give a different presentation
of the Tetrahedron algebra by generators and relations, inspired in
the properties of the $\calA$-basis used throughout. Finally,
Section 6 will deal with the universal central extension of the
Tetrahedron algebra studied in \cite{GeorgiaPaul}.

\bigskip

\section{The Tetrahedron algebra and the three point $\frsldos$ loop
algebra}

The Tetrahedron algebra $\frgtetra$ has been defined in \cite{Paul}.
It is the Lie algebra over $k$ with generators
\begin{equation}\label{eq:Xij}
\{ X_{ij}: i,j\in\{0,1,2,3\},\, i\ne j\}
\end{equation}
and the relations
\begin{subequations}\label{eq:rel}
\begin{align}
&X_{ij}+X_{ji}=0\quad\text{for $i\ne j$,}\label{eq:rel1}\\
&[X_{ij},X_{jk}]=2(X_{ij}+X_{jk})\quad\text{for mutually distinct
$i,j,k$,}\label{eq:rel2}\\
&[X_{hi},[X_{hi},[X_{hi},X_{jk}]]]=4[X_{hi},X_{jk}]\quad \text{for
mutually distinct $h,i,j,k$.}\label{eq:rel3}
\end{align}
\end{subequations}

One of the main results in \cite{Paul} relates the Tetrahedron
algebra to the three point $\frsldos$ loop algebra
$\frg=\frsldos\otimes \calA$, where $\frsldos$ is the Lie algebra of
two by two traceless matrices over $k$ and $\calA$ is the unital
commutative associative algebra $k[t,t^{-1},(1-t)^{-1}]$ ($t$ an
indeterminate), which is a subalgebra of the field of fractions
$k(t)$. To present the precise relationship in \cite[Proposition 6.5
and Theorem 11.5]{Paul}, first consider the basis
$\bigl\{x=\bigl(\begin{smallmatrix}
-1&2\\0&1\end{smallmatrix}\bigr),\,
 y=\bigl(\begin{smallmatrix}
-1&0\\-2&1\end{smallmatrix}\bigr),\,
 z=\bigl(\begin{smallmatrix}
1&0\\0&-1\end{smallmatrix}\bigr)\bigr\}$ of $\frsldos$, whose
elements satisfy $[x,y]=2(x+y)$, $[y,z]=2(y+z)$, $[z,x]=2(z+x)$.
Then, because of its own definition by generators and relations,
there is a Lie algebra homomorphism \cite[Proposition 6.5]{Paul}
\[
\Psi:\frgtetra\longrightarrow \frg=\frsldos\otimes\calA
\]
determined by:
\begin{equation}\label{eq:Psi}
\begin{split}
\Psi(X_{12})=&x\otimes 1,\quad \Psi(X_{23})=y\otimes 1,\quad
   \Psi(X_{31})=z\otimes 1,\\
&\Psi(X_{03})=y\otimes t +z\otimes (t-1),\\
&\Psi(X_{01})=z\otimes t'+x\otimes (t'-1),\\
&\Psi(X_{02})=x\otimes t''+y\otimes (t''-1),
\end{split}
\end{equation}
where $t'=1-t^{-1}$ and $t''=(1-t)^{-1}$ (see \cite[Lemma
6.2]{Paul}). This homomorphism $\Psi$ was proved to be an
isomorphism in \cite[Theorem 11.5]{Paul}. Here this will be proved
in another way. Some useful results will come out during this
process.

The symmetric group
\[
S_4=\{\sigma:\{0,1,2,3\}\rightarrow\{0,1,2,3\}: \sigma\text{\ is a
bijective map}\}
\]
embeds naturally in the group of automorphisms $\Aut(\frgtetra)$ by
means of
\[
\sigma(X_{ij})=X_{\sigma(i)\sigma(j)},
\]
for any $\sigma\in S_4$ and $0\leq i\ne j\leq 3$. (Here the actions
will always be taken on the left.)

Consider the following generators of $S_4$:
\[
\begin{split}
&\tau_1=(12)(30): 1\leftrightarrow 2,\ 3\leftrightarrow 0,\\
&\tau_2=(23)(10): 2\leftrightarrow 3,\ 1\leftrightarrow 0,\\
&\varphi=(123): 1\mapsto 2\mapsto 3\mapsto 1,\ 0\mapsto 0,\\
&\tau=(12): 1\leftrightarrow 2,\ 0\mapsto 0,\ 3\mapsto 3.
\end{split}
\]
The elements $\tau_1$ and $\tau_2$ generate Klein's $4$ group, while
$\varphi$ and $\tau$ generate a copy of the symmetric group $S_3$
(recall that $S_4$ is the semidirect product of these two
subgroups).

\begin{theorem}\label{th:automorphisms}
$S_4$ embeds as a group of automorphisms of
$\frg=\frsldos\otimes\calA$ in the following way:
\begin{romanenumerate}
\item $\varphi=\varphi_\frs\otimes\varphi_\calA$, where
$\varphi_\frs$ is the order $3$ automorphism of $\frsldos$ given by
\[
\varphi_\frs(x)=y,\quad \varphi_\frs(y)=z,\quad \varphi_\frs(z)=x,
\]
and $\varphi_\calA$ is the order $3$ automorphism of the $k$-algebra
$\calA$ determined by
\[
\varphi_\calA(t)=1-t^{-1}=t'.
\]
In particular, $\varphi$ is an $\calA$-semilinear automorphism of
$\frg$ and $\varphi_\calA$ is its associated automorphism of
$\calA$. That is, $\varphi(ga)=\varphi(g)\varphi_\calA(a)$ for any
$g\in \frg$ and $a\in\calA$.

\item $\tau=\tau_\frs\otimes\tau_\calA$, where $\tau_\frs$ is the
order $2$ automorphism of $\frsldos$ given by
\[
\tau_\frs(x)=-x,\quad \tau_\frs(y)=-z,\quad \tau_\frs(z)=-y,
\]
and $\tau_\calA$ is the order $2$ automorphism of $\calA$ determined
by $\tau_\calA(t)=1-t$. In particular, $\tau$ is an
$\calA$-semilinear automorphism of $\frg$ and $\tau_\calA$ is its
associated automorphism of $\calA$.

\item $\tau_1$ is the automorphism of $\frg$, as a Lie algebra over
$\calA$, given by
\begin{equation}\label{eq:tau1}
\begin{split}
\tau_1(x\otimes 1)&=-x\otimes 1,\\
\tau_1(y\otimes 1)&=-\bigl(z\otimes t'+x\otimes (t'-1)\bigr),\\
\tau_1(z\otimes 1)&=x\otimes t''+y\otimes (t''-1).
\end{split}
\end{equation}

\item $\tau_2$ is the automorphism of $\frg$, as a Lie algebra over
$\calA$, given by
\begin{equation}\label{eq:tau2}
\begin{split}
\tau_2(x\otimes 1)&=y\otimes t+z\otimes (t-1),\\
\tau_2(y\otimes 1)&=-y\otimes 1,\\
\tau_2(z\otimes 1)&=-\bigl(x\otimes t''+y\otimes (t''-1)\bigr).
\end{split}
\end{equation}
\end{romanenumerate}
Moreover, under these actions of $S_4$ on $\frgtetra$ and on $\frg$,
the homomorphism $\Psi$ in \eqref{eq:Psi} becomes a homomorphism of
Lie algebras with $S_4$-action. That is,
\begin{equation}\label{eq:PsisigmaX}
\Psi\bigl(\sigma(X)\bigr)=\sigma\bigl(\Psi(X)\bigr),
\end{equation}
for any $\sigma\in S_4$ and $X\in\frgtetra$.
\end{theorem}

\begin{proof}
It is easy to check that $\varphi_\frs$ and $\tau_\frs$ are
automorphisms of $\frsldos$ of order $3$ and $2$, respectively, and
that $\varphi_\frs\tau_\frs=\tau_\frs\varphi_\frs^2$. Also,
$\varphi_\calA\tau_\calA=\tau_\calA\varphi_\calA^2$. Hence $S_3$
embeds in $\Aut\frg$ by identifying $\varphi$ to
$\varphi_\frs\otimes\varphi_\calA$ and $\tau$ to
$\tau_\frs\otimes\tau_\calA$.

Now, equation \eqref{eq:tau1} defines a unique Lie algebra
homomorphism of $\frg$, as a Lie algebra over $\calA$, which
satisfies $\tau_1^2=1$. The same happens for $\tau_2$ and
$\tau_1\tau_2=\tau_2\tau_1$. For instance, using that
$(t'-1)t=-1=t''(t-1)$, we check
\[
\begin{split}
\tau_1\tau_2(x\otimes 1)&=
  \tau_1\bigl(y\otimes t+ z\otimes (t-1)\bigr)\\
  &=-z\otimes t't-x\otimes (t'-1)t+x\otimes t''(t-1)+y\otimes
  (t''-1)(t-1)\\
  &=-z\otimes (t-1)+x\otimes 1-x\otimes 1 -y\otimes t\\
  &=-\bigl(y\otimes t+z\otimes(t-1)\bigr)\\
  &=-\tau_2(x\otimes 1)=\tau_2\tau_1(x\otimes 1).
\end{split}
\]
Now, it has to be checked that $\varphi\tau_1=\tau_2\varphi$,
$\varphi\tau_2=\tau_1\tau_2\varphi$, $\tau\tau_1=\tau_1\tau$ and
$\tau\tau_2=\tau_1\tau_2\tau$. In all cases, the maps on both sides
are $\calA$-semilinear with the same associated automorphisms of
$\calA$, so it is enough to check the equalities on the basis
$\{x\otimes 1,y\otimes 1,z\otimes 1\}$ of $\frg$ over $\calA$. This
is straightforward.

For the last part, it is enough to check \eqref{eq:PsisigmaX} for
$X=X_{ij}$, $i\ne j$ in $\{0,1,2,3\}$, and again this is a routine
verification.
\end{proof}

\medskip

In order to show that $\Psi$ is an isomorphism, it is better to work
with a different basis of $\frg$ over $\calA$. Consider the
following elements of $\frg$:
\begin{equation}\label{eq:Abasis}
\begin{split}
u_0&=\frac{1}{4}\Psi(X_{02}+X_{31})
  =\frac{1}{4}\bigl(z\otimes 1+x\otimes t''+y\otimes(t''-1)\bigr),\\
u_1&=\frac{1}{4}\Psi(X_{03}+X_{12})
  =\frac{1}{4}\bigl(x\otimes 1+y\otimes t+z\otimes(t-1)\bigr),\\
u_2&=\frac{1}{4}\Psi(X_{01}+X_{23})
  =\frac{1}{4}\bigl(y\otimes 1+z\otimes t'+x\otimes(t'-1)\bigr).
\end{split}
\end{equation}
Observe that, since $\varphi_\calA(t)=t'$ and
$\varphi_\calA(t')=t''$, these elements are permuted cyclically by
the order $3$ automorphism $\varphi$ in Theorem
\ref{th:automorphisms}.

\begin{theorem}\label{th:Abasis}
With $u_0$, $u_1$ and $u_2$ as above:
\begin{romanenumerate}
\item $\{u_0,u_1,u_2\}$ is a basis of $\frg$ as a module over
$\calA$.

\item $[u_0,u_1]=-u_2t$, $[u_1,u_2]=-u_0t'$, $[u_2,u_0]=-u_1t''$.

\item $u_0$, $u_1$ and $u_2$ generate $\frg$ as a Lie algebra over
$k$.
\end{romanenumerate}
\end{theorem}
\begin{proof}
Let us start with item (ii). Use $\Psi$ and the relation
\eqref{eq:rel2} to get
\[
\begin{split}
[u_0,u_1]&=\frac{1}{16}\Psi\Bigl(
           [X_{02}+X_{31},X_{03}+X_{12}]\Bigr)\\
 &=\frac{1}{16}\Psi\Bigl(-[X_{20},X_{03}]-[X_{02},X_{21}]
        +[X_{13},X_{30}]+[X_{31},X_{12}]\Bigr)\\
 &=\frac{1}{16}\Psi\Bigl(-2(X_{20}+X_{03})-2(X_{02}+X_{21})\\
 &\qquad\qquad\qquad
        +2(X_{13}+X_{30})+2(X_{31}+X_{12})\Bigr)\\
 &=\frac{1}{16}\Psi\Bigl(-4X_{03}+4X_{12}\Bigr)
    =-\frac{1}{4}\Psi\bigl(X_{03}-X_{12}\bigr)\\
 &=-\frac{1}{4}\Psi\bigl(-x\otimes 1+y\otimes t+z\otimes(t-1)\bigr).
\end{split}
\]
On the other hand,
\[
\begin{split}
-u_2t&=-\frac{1}{4}\bigl(y\otimes t+z\otimes t't+x\otimes
           (t'-1)t\bigr)\\
 &=-\frac{1}{4}\bigl(y\otimes t+z\otimes (t-1)-x\otimes 1\bigr),
\end{split}
\]
since $(t'-1)t=-1$. Hence $[u_0,u_1]=-u_2t$. Now apply $\varphi$ to
get $[u_1,u_2]=[\varphi(u_0),\varphi(u_1)]=-\varphi(u_2t)
=-\varphi(u_2)\varphi_\calA(t)=-u_0t'$ and, in the same vein,
$[u_2,u_0]=-u_1t''$, as required.

To prove (i), first note that $\frg=\frsldos\otimes\calA$ is a
subalgebra of the simple Lie algebra $\frsldos\otimes k(t)$ over the
field $k(t)$. The $k(t)$-subalgebra $\frs$ generated by $u_0$, $u_1$
and $u_2$ is perfect ($\frs=[\frs,\frs]$) by (ii), and hence its
dimension cannot be $\leq 2$ (otherwise, it would be a solvable Lie
algebra). Therefore, $\{u_0,u_1,u_2\}$ is a basis of
$\frsldos\otimes k(t)$ over $k(t)$. In particular, $u_0$, $u_1$ and
$u_2$ are linearly independent over $\calA$, and $u_0\calA\oplus
u_1\calA\oplus u_2\calA$ is an $\calA$-subalgebra of $\frg$.
Moreover, the computations above show that
\begin{equation}\label{eq:u2t}
4u_2t=y\otimes t+z\otimes(t-1)-x\otimes 1,
\end{equation}
while $4u_1=x\otimes 1 + y\otimes t+z\otimes(t-1)$. Hence,
\begin{equation}\label{eq:xotimes1}
x\otimes 1=2(u_1-u_2t)
\end{equation}
belongs to $u_0\calA\oplus u_1\calA\oplus u_2\calA$. Apply $\varphi$
to get $y\otimes 1=2(u_2-u_0t')\in u_0\calA\oplus u_1\calA\oplus
u_2\calA$ and also $z\otimes 1=2(u_0-u_1t'')\in u_0\calA\oplus
u_1\calA\oplus u_2\calA$. This proves (i).

Finally, let us denote now by $\frs$ the $k$-subalgebra of $\frg$
generated by $u_0$, $u_1$ and $u_2$. Observe that $\frs$ is
invariant under the action of the order $3$ automorphism $\varphi$,
so it is enough to prove that $u_0\calA$ is contained in $\frs$.
From (ii) we obtain
\begin{multline*}
[u_1,[u_1,u_0t^n]]=[u_1,u_2t^{n+1}]=-u_0t't^{n+1}\\
=-u_0t^n(t-1)=-u_0t^{n+1}+u_0t^n
\end{multline*}
since $t't=t-1$. Hence, an induction argument shows that $u_0t^n$ is
in $\frs$ for any $n\geq 0$. In the same vein,
\[
[u_2,[u_2,u_0(t')^n]]=[u_2,u_1(t')^nt'']=-u_0(t')^{n+1}t''=-u_0(t')^n(t'-1),
\]
and $u_0(t')^n\in\frs$ for any $n\geq 0$. Finally, since $u_0t^n$ is
in $\frs$ for any $n\geq 0$, and $\frs$ is invariant under
$\varphi$, we get $u_2(t'')^n\in\frs$ too, for any $n\geq 0$. Hence:
\[
[u_1,u_2(t'')^n(1-t'')]=-u_0(t'')^n(1-t'')t'=-u_0(t'')^n\in\frs
\]
for any $n\geq 0$. But $\{1\}\cup\{ t^n,(t')^n,(t'')^n : n\in\bN\}$
is a basis of $\calA$ over $k$ (\cite[Lema 6.3]{Paul}), so
$u_0\calA$ is contained in $\frs$, as required.
\end{proof}

\smallskip

\begin{corollary}\label{co:Psionto}
The homomorphism $\Psi$ is onto.
\end{corollary}

\smallskip

As in \cite{Paul}, let $\Omega$ (respectively $\Omega'$, $\Omega''$)
denote the subalgebra of $\frgtetra$ generated by $X_{12}$ and
$X_{03}$ (respectively $X_{23}$ and $X_{01}$, $X_{31}$ and
$X_{02}$). Note that $\Omega'=\varphi(\Omega)$ and
$\Omega''=\varphi(\Omega')$. In \cite[Proposition 7.8]{Paul} it is
proved that $\frgtetra$ is the direct sum of the subalgebras
$\Omega$, $\Omega'$ and $\Omega''$. A simpler proof can be given as
follows:

\begin{lemma}\label{le:previous}
Let $S_1$ and $S_2$ be two subspaces of a Lie algebra $\frl$ such
that $[S_1,S_2]\subseteq S_1+S_2$ holds, and let $\frs_i$ be the
subalgebra generated by $S_i$, $i=1,2$. Then
$[\frs_1,\frs_2]\subseteq \frs_1+\frs_2$. In particular,
$\frs_1+\frs_2$ is a subalgebra of $\frl$.
\end{lemma}
\begin{proof}
 From $[S_1,S_2]\subseteq S_1+S_2$, it follows that
$[S_1,\frs_2]\subseteq S_1+\frs_2$, and then that
$[\frs_1,\frs_2]\subseteq \frs_1+\frs_2$.
\end{proof}

\begin{proposition}\label{pr:Omegaprimereprime}
$\frgtetra=\Omega + \Omega' +\Omega''$.
\end{proposition}
\begin{proof}
Let $S$ (respectively $S'$, $S''$) denote the subspace spanned by
$X_{12}$ and $X_{03}$ (respectively $X_{23}$ and $X_{01}$, $X_{31}$
and $X_{02}$). Then, by \eqref{eq:rel2}, $[S,S']\subseteq S+S'$, so
the previous Lemma gives $[\Omega,\Omega']\subseteq \Omega+\Omega'$,
and, similarly, $[\Omega',\Omega'']\subseteq \Omega'+\Omega''$ and
$[\Omega'',\Omega]\subseteq \Omega''+\Omega$. Therefore,
$\Omega+\Omega'+\Omega''$ is a subalgebra of $\frgtetra$, which
contains all the generators $X_{ij}$, so it is the whole
$\frgtetra$.
\end{proof}

\smallskip

The images under $\Psi$ of these subalgebras are given in the next
result:

\begin{proposition}\label{pr:Omega}
\null\quad\null
\begin{romanenumerate}
\item $\Psi(\Omega)=u_0(t-1)k[t]\oplus u_1k[t]\oplus u_2tk[t]$,
\item $\Psi(\Omega')=u_0t'k[t']\oplus u_1(t'-1)k[t']\oplus u_2k[t']$,
\item $\Psi(\Omega'')=u_0k[t'']\oplus u_1t''k[t'']\oplus u_2(t''-1)k[t'']$.
\end{romanenumerate}
In particular, $\frg$ is the direct sum of the subalgebras
$\Psi(\Omega)$, $\Psi(\Omega')$, and $\Psi(\Omega'')$.
\end{proposition}
\begin{proof}
(ii) and (iii) are obtained from (i) by applying the order $3$
automorphism $\varphi$. To prove (i), first note that since
$tt'=t-1$ and $t''(t-1)=-1$, the $k$-subspace $u_0(t-1)k[t]\oplus
u_1k[t]\oplus u_2tk[t]$ is a $k$-subalgebra of $\frg$, and it
contains $u_1=\frac{1}{4}\Psi(X_{03}+X_{12})$ and
$u_2t=\frac{1}{4}\Psi(X_{03}-X_{12})$ \eqref{eq:u2t}, so it
certainly contains $\Psi(\Omega)$. But one has
$\calA=(t-1)k[t]\oplus t'k[t']\oplus k[t'']$ because of \cite[Lemma
6.3]{Paul}, so applying $\varphi_\calA$, also $\calA=k[t]\oplus
(t'-1)k[t']\oplus t''k[t'']= tk[t]\oplus k[t']\oplus (t''-1)k[t'']$.
Hence $\frg$ is the direct sum of the subalgebras on the right hand
sides of (i), (ii), and (iii). Since $\Psi$ is onto (Corollary
\ref{co:Psionto}) and $\frgtetra=\Omega+\Omega'+\Omega''$
(Proposition \ref{pr:Omegaprimereprime}), the result follows.
\end{proof}

\smallskip

The next result simplifies the work to prove that $\Psi$ is
one-to-one:

\begin{lemma}\label{le:PsiOmega}
If the restriction of $\Psi$ to $\Omega$ is one-to-one, so is
$\Psi$.
\end{lemma}
\begin{proof}
Assume that the restriction of $\Psi$ to $\Omega$ is one-to-one.
Then, since $\Psi\bigl(\Omega\cap(\Omega'+\Omega'')\bigr)$ is
contained in $\Psi(\Omega)\cap(\Psi(\Omega')+\Psi(\Omega''))=0$ (as
$\frg$ is the direct sum of $\Psi(\Omega)$, $\Psi(\Omega')$, and
$\Psi(\Omega'')$ by Proposition \ref{pr:Omega}), it follows that
$\Omega\cap(\Omega'+\Omega'')$ is contained in the kernel of the
restriction $\Psi\vert_{\Omega}$, which is assumed to be $0$. Hence,
by the cyclic symmetry provided by $\varphi$,
$\frgtetra=\Omega\oplus\Omega'\oplus\Omega''$ and, again using the
cyclic symmetry, the restriction of $\Psi$ to each of these three
direct summands is one-to-one, and so is $\Psi$.
\end{proof}

\smallskip

Therefore, it is enough to prove that
$\Psi\vert_\Omega:\Omega\rightarrow \frg$ is one-to-one. But
$\Omega$ is the homomorphic image of the Onsager algebra $\calO$
(see \cite[Section 4]{Paul}), which is the Lie algebra over $k$ with
generators $A$ and $B$ and relations
\[
[A,[A,[A,B]]]=4[A,B],\quad [B,[B,[B,A]]]=4[B,A],
\]
under the homomorphism which takes $A$ to $X_{12}$ and $B$ to
$X_{03}$.

Hence, in order to prove that $\Psi\vert_\Omega$ is one-to-one, it
is enough to prove the following:

\begin{lemma}\label{le:phiO}
The Lie algebra homomorphism $\phi:\calO\rightarrow \frg$ determined
by $\phi(A)=\Psi(X_{12})$ and $\phi(B)=\Psi(X_{03})$ is one-to-one.
\end{lemma}
\begin{proof}
First note that $\phi(A)=\Psi(X_{12})=x\otimes 1=2(u_1-u_2t)$
\eqref{eq:xotimes1}, and $\phi(B)=\Psi(X_{03})=y\otimes t+z\otimes
(t-1)=2(u_1+u_2t)$ \eqref{eq:u2t}. Besides \cite{Onsager}, $\calO$
has a basis $\{ A_m: m\in\bZ\}\cup\{ G_l: l\in\bN\}$, where $A_0=A$,
$A_1=B$ and
\[
[A_l,A_m]=2G_{l-m}\ (l>m),\quad [G_l,A_m]=A_{m+l}-A_{m-l},\quad
[G_l,G_m]=0.
\]
Denote by $\calO_A$ (respectively $\calO_G$) the linear span of
$\{A_m:m\in\bZ\}$ (respectively $\{G_l: l\in\bN\}$). Then
$\ad_{G_1}\vert_{\calO_A}:\calO_A\rightarrow \calO_A$ is one-to-one,
and so is $\ad_{A_0}\vert_{\calO_G}:\calO_G\rightarrow \calO_A$.
This shows that any nonzero ideal of $\calO$ intersects nontrivially
$\calO_A$. Also, $\{(\ad_{G_1})^m(A_0),(\ad_{G_1})^m(A_1): m\geq
0\}$ is a basis of $\calO_A$, and so is
$\{(\ad_{G_1})^m(A_0+A_1),(\ad_{G_1})^m(A_0-A_1):m\geq 0\}$.

Hence, in order to prove that $\phi$ is one-to-one, it is enough to
prove that so is $\phi\vert_{\calO_A}$ and, hence, to prove that
$\{\phi\bigl((\ad_{G_1})^m(A_0+A_1)\bigr),\phi\bigl((\ad_{G_1})^m(A_0-A_1)\bigr):
m\geq 0\}$ is a linearly independent set in $\frg$. But,
\[
\begin{split}
\phi(G_1)&=\frac{1}{2}\phi([A_1,A_0])
    =\frac{1}{2}[\Psi(X_{03}),\Psi(X_{12})]\\
    &=2[u_1+u_2t,u_1-u_2t]=-4[u_1,u_2t]=4u_0tt'=4u_0(t-1),\\[6pt]
\phi(A_0+A_1)&=2(u_1-u_2t)+2(u_1+u_2t)=4u_1,\\
\phi(A_0-A_1)&=2(u_1+u_2t)-2(u_1-u_2t)=4u_2t,
\end{split}
\]
so we must check that
$\{(\ad_{u_0(t-1)})^m(u_1),(\ad_{u_0(t-1)})^m(u_2t):m\geq 0\}$ is a
linearly independent set. Now,
\[
\begin{split}
[u_0(t-1),u_1]&=-u_2(t-1),\\
[u_0(t-1),u_2t]&=u_1t''(t-1)t=-u_1t,
\end{split}
\]
so
\[
\begin{split}
(\ad_{u_0(t-1)})^{2m}(u_1)&=u_1\bigl(t(t-1)\bigr)^m,\\
(\ad_{u_0(t-1)})^{2m}(u_2t)&=(u_2t)\bigl(t(t-1)\bigr)^m,\\
(\ad_{u_0(t-1)})^{2m+1}(u_1)&=-(u_2t)(t-1)\bigl(t(t-1)\bigr)^m,\\
(\ad_{u_0(t-1)})^{2m+1}(u_2t)&=-u_1t\bigl(t(t-1)\bigr)^m,
\end{split}
\]
and all these elements are linearly independent over $k$.
\end{proof}


\begin{corollary}\label{co:Main}
(See \cite[Theorems 11.5 and 11.6, and Corollary 12.5]{Paul}.)
\begin{romanenumerate}
\item $\Psi:\frgtetra\rightarrow \frg$ is an isomorphism.
\item $\frgtetra$ is the direct sum of its subalgebras $\Omega$,
$\Omega'$ and $\Omega''$.
\item $\Omega$ is isomorphic to the Onsager algebra.
\end{romanenumerate}
\end{corollary}
\begin{proof}
(i) follows from the results above. Then (ii) follows from
Proposition \ref{pr:Omega}, and (iii) follows because the
epimorphism $\calO\rightarrow \Omega$ such that $A\mapsto X_{12}$
and $B\mapsto X_{03}$ is one-to-one by Lemma \ref{le:phiO}.
\end{proof}

\bigskip

\section{The solution to \cite[Problem 13.4]{Paul}}

The action of $S_4$ on $\frgtetra$ (or on
$\frg=\frsldos\otimes\calA$) restricts to an action of Klein's $4$
group, which gives a $\bZ_2\times\bZ_2$-grading on $\frg$, as in
\cite[(1.1)]{EO}:
\begin{equation}\label{eq:tg0g1g2}
\frg=\frt\oplus\frg_0\oplus\frg_1\oplus\frg_2,
\end{equation}
where
\[
\begin{split}
\frt&=\{g\in\frg : \tau_1(g)=g,\, \tau_2(g)=g\}\
       (=\frg_{(\bar 0,\bar 0)}),\\
\frg_0&=\{g\in\frg : \tau_1(g)=g,\, \tau_2(g)=-g\}\
        (=\frg_{(\bar 1,\bar 0)}),\\
\frg_1&=\{g\in\frg : \tau_1(g)=-g,\, \tau_2(g)=g\}\
        (=\frg_{(\bar 0,\bar 1)}),\\
\frg_2&=\{g\in\frg : \tau_1(g)=-g,\, \tau_2(g)=-g\}\
        (=\frg_{(\bar 1,\bar 1)}).
\end{split}
\]

In \cite[Problem 13.4]{Paul} it is posed the question of proving
that $\frt=0$, of obtaining a basis for each of these subspaces
$\frg_i$, and of investigating the relationship between this
decomposition \eqref{eq:tg0g1g2} and the decomposition
$\frgtetra=\Omega\oplus\Omega'\oplus\Omega''$ in Corollary
\ref{co:Main}.

The use of the $\calA$-basis $\{u_0,u_1,u_2\}$ of $\frg$ (Theorem
\ref{th:Abasis}) makes the determination of the subspaces in
\eqref{eq:tg0g1g2} quite easy:

\begin{theorem}\label{th:tg0g1g2} With the previous notations,
\[
\frt=0,\quad \frg_0=u_0\calA,\quad \frg_1=u_1\calA,\quad
\frg_2=u_2\calA.
\]
\end{theorem}
\begin{proof}
Recall that $\tau_1=(12)(03)$ and $\tau_2=(23)(10)$. Since
$\tau_1(X_{02})=X_{31}$ and $\tau_2(X_{02})=X_{13}=-X_{31}$, it
follows that $u_0=\frac{1}{4}\Psi(X_{02}+X_{31})$ belongs to
$\frg_0$. But the automorphisms $\tau_1$ and $\tau_2$ are
$\calA$-linear (Proposition \ref{th:automorphisms}), so $u_0\calA$
is contained in $\frg_0$. In the same vein it is proved that
$u_1\calA\subseteq \frg_1$ and $u_2\calA\subseteq \frg_2$; and
Theorem \ref{th:Abasis} finishes the proof.
\end{proof}

Since the set $\{1\}\cup \{t^n,(t')^n,(t'')^n: n\in\bN\}$ is a
$k$-basis of $\calA$ (\cite[Lemma 6.3]{Paul}), the following result,
which solves part of \cite[Problem 13.4]{Paul}, is clear:

\begin{corollary}\label{co:basis}
For $i=0,1,2$, the set $\{u_i\}\cup\{ u_it^n,u_i(t')^n,u_i(t'')^n:
n\in \bN\}$ is a $k$-basis of the space $\frg_i$.
\end{corollary}

\smallskip

Actually, \cite[Lemma 6.3]{Paul} shows that $\calA=(t-1)k[t]\oplus
t'k[t']\oplus k[t'']$, so
\begin{multline*}
\frg_0=u_0(t-1)k[t]\oplus u_0t'k[t']\oplus u_0k[t'']\\
=
 \bigl(\Psi(\Omega)\cap \frg_0\bigr)\oplus
 \bigl(\Psi(\Omega')\cap \frg_0\bigr)\oplus
 \bigl(\Psi(\Omega'')\cap \frg_0\bigr),
\end{multline*}
(see Proposition \ref{pr:Omega}), and something similar holds for
$\frg_1$ and $\frg_2$. This gives the relationship between the
decompositions $\frg=\frg_0\oplus\frg_1\oplus\frg_2$ and the
decomposition $\frgtetra=\Omega\oplus\Omega'\oplus\Omega''$.

\smallskip

\begin{remark}
Also, since $t''(t-1)=-1$ and $(t''-1)(t-1)=-t$, it follows that
\[
\begin{split}
u_0(t-1)k[t]&=\Bigl(z\otimes(t-1)+x\otimes t''(t-1)+y\otimes
(t''-1)(t-1)\Bigr) k[t]\\
 &\subseteq x\otimes k[t]\oplus y\otimes
tk[t]\oplus z\otimes (t-1)k[t].
\end{split}
\]
Besides, the definition of $u_1$ shows immediately that
\[
u_1k[t]\subseteq x\otimes k[t]\oplus y\otimes tk[t]\oplus z\otimes
(t-1)k[t],
\]
and since $t't=t-1$ and $(t'-1)t=-1$, also
\[
u_2tk[t]\subseteq x\otimes k[t]\oplus y\otimes tk[t]\oplus z\otimes
(t-1)k[t].
\]
In this way, one recovers \cite[Corollary 13]{Paul}, which asserts
that
\[
\Psi(\Omega)= x\otimes k[t]\oplus y\otimes tk[t]\oplus z\otimes
(t-1)k[t].
\]
\end{remark}

\bigskip

\section{The normal Lie related triple algebra associated to the
Tetrahedron algebra}

Following \cite{EO}, given the Lie algebra $\frg$ on which $S_4$
acts as automorphisms, there exists a structure of normal Lie
related triple algebra defined on $\frg_0$, which essentially
determines $\frg$.

A \emph{normal Lie related triple algebra} $(A,\cdot,\bar{\ })$ is
an algebra with multiplication $\cdot$, with an involution $\bar{\
}$, and endowed with a skew-symmetric bilinear map $ \delta:A\times
A\rightarrow \lrt(A,\cdot,\bar{\ }) $, where
\begin{multline*}
\lrt(A,\cdot,\bar{\ })=\{(d_0,d_1,d_2)\in \frgl(A)^3 : \bar
d_i(x\cdot y)
  =d_{i+1}(x)\cdot y+x\cdot d_{i+2}(y)\\
  \text{for any }x,y\in A\text{ and }i\in\bZ_3\},
\end{multline*}
satisfying some conditions (see \cite[(2.34)]{Okubo} for a complete
definition).

Here $\frg_0=u_0\calA$ can be identified with $\calA$ by means of
$\iota_0:\calA\rightarrow \frg_0$, $a\mapsto \iota_0(a)=u_0a$. Then,
according to \cite{EO}, one has to consider the identifications
$\iota_i:\calA\rightarrow \frg_i$ given by
\[
\iota_i(a)=\varphi^{i}\bigl(\iota_0(a)\bigr)=\varphi^i(u_0a)
 =\varphi^i(u_0)\varphi_\calA^i(a)=u_i\varphi_\calA^i(a),
\]
for $i=0,1,2$ (see Proposition \ref{th:automorphisms}). Therefore,
by Theorem \ref{th:tg0g1g2},
\[
\frg=\oplus_{i=0}^2\iota_i(\calA),
\]
and the action of $\calA$ on $\frg$ is given by
\[
\iota_i(a)b=u_i\varphi_\calA^i(a)b
 =u_i\varphi_\calA^i\bigl(a\varphi_\calA^{-i}(b)\bigr)=
 \iota_i\bigl(a\varphi_\calA^{-i}(b)\bigr),
\]
for any $a,b\in\calA$ and $i=0,1,2$.

Now, according to \cite[Section 2]{EO}, the $k$-vector space $\calA$
is endowed with an involution $a\mapsto \bar a$ determined by
$\iota_0(\bar a)=-\tau\bigl(\iota_0(a)\bigr)$. That is, by
Proposition \ref{th:automorphisms},
\[
\begin{split}
u_0\bar a&=-\tau(u_0a)=-\tau(u_0)\tau_\calA(a)\\
 &=-\frac{1}{4}\tau\bigl(z\otimes 1+x\otimes t''+y\otimes
 (t''-1)\bigr)\tau_\calA(a)\\
 &=-\frac{1}{4}\bigl(\tau_\frs(z)\otimes 1+\tau_\frs(x)\otimes
 \tau_\calA(t'')
 +\tau_\frs(y)\otimes\tau_\calA(t''-1)\bigr)\tau_\calA(a)\\
 &=-\frac{1}{4}\bigl(-y\otimes 1-x\otimes(1-t')+z\otimes
 t'\bigr)\tau_\calA(a)\\
 &=-u_0t'\tau_\calA(a),
\end{split}
\]
as $t''t'=t'-1$. Therefore,
\[
\bar a=-t'\tau_\calA(a)
\]
for any $a\in \calA$. Note that since $\tau^2=1$, $\bar{\bar a}=a$
for any $a\in \calA$.

Also, $\calA$ is endowed with a multiplication determined by
\[
\iota_0(\overline{a\cdot b})=[\iota_1(a),\iota_2(b)].
\]
 Hence,
\[
\begin{split}
u_0(\overline{a\cdot b})&=
 \iota_0(\overline{a\cdot b})=[\iota_1(a),\iota_2(b)]\\
 &=[u_1\varphi_\calA(a),u_2\varphi_\calA^2(b)]\\
 &=[u_1,u_2]\varphi_\calA(a)\varphi_\calA^2(b)\\
 &=-u_0t'\varphi_\calA(a)\varphi_\calA^2(b)
\end{split}
\]
(see Theorem \ref{th:Abasis}), so
\[
\begin{split}
a\cdot b&=\overline{-t'\varphi_\calA(a)\varphi_\calA^2(b)}\\
 &=t'\tau_\calA\bigl(t'\varphi_\calA(a)\varphi_\calA^2(b)\bigr)\\
 &=\bigl(\tau_\calA\varphi_\calA(a)\bigr)
    \bigl(\tau_\calA\varphi_\calA^2(b)\bigr),
\end{split}
\]
for any $a,b\in\calA$ (since $\tau_\calA(t')=(t')^{-1}$).

Summarizing, we have obtained:

\begin{proposition} The normal Lie related triple algebra associated
to the $S_4$-action on the Lie algebra $\frg=\frsldos\otimes\calA$
is isomorphic to $(\calA,\cdot, \bar\ )$, where
\[
\left\{\begin{aligned} &a\cdot b
         =\bigl(\tau_\calA\varphi_\calA(a)\bigr)
             \bigl(\tau_\calA\varphi_\calA^2(b)\bigr)\\
         &\bar a=-t'\tau_\calA(a)\end{aligned}\right.
\]
for any $a,b\in\calA$.
\end{proposition}

\medskip

This result highlights a family of normal Lie related triple
algebras, whose associated Lie algebras with $S_4$-action satisfy
that there are no nonzero elements fixed by Klein's $4$ group:

\begin{proposition} Let $A$ be a unital commutative associative
algebra endowed with a group homomorphism $\theta:S_3\rightarrow
\Aut(A)$ (the image of any $\sigma$ under $\theta$ will be denoted
by $\sigma$ too), and an element $s\in A$ such that
$s\varphi(s)\varphi^2(s)=1=s\tau(s)$ (notation as in Section 1).
Define on $A$ a new multiplication by
\[
a\cdot b=\bigl(\tau\varphi(a)\bigr)\bigl(\tau\varphi^2(b)\bigr)
\]
for any $a,b\in A$, and a linear map $A\rightarrow A$, $a\mapsto
\bar a$, by
\[
\bar a=s\tau(a).
\]
Then $(A,\cdot,\bar\ )$ is a normal Lie related triple algebra with
trivial associated bilinear map $\delta: A\times A\rightarrow
\lrt(A,\cdot,\bar\ )$.
\end{proposition}
\begin{proof}
Let us first check that $a\mapsto \bar a$ is an involution of
$(A,\cdot)$. For any $a\in A$,
\[
\bar{\bar a}=s\tau(\bar a)=s\tau(s\tau(a))=s\tau(s)\tau^2(a)=a,
\]
while for any $a,b\in A$,
\[
\begin{split}
\bar b\cdot \bar a&=
 \bigl(\tau\varphi(\bar b)\bigr)\bigl(\tau\varphi^2(\bar a)\bigr)\\
 &=\bigl(\tau\varphi(s\tau(b))\bigr)\bigl(\tau\varphi^2(s\tau(a))\bigr)\\
 &=\tau\bigl(\varphi(s)\varphi^2(s)\bigr)\varphi^2(b)\varphi(a)\quad
 \text{(as $\tau\varphi\tau=\varphi^2$)}\\
 &=s\tau\bigl(s\varphi(s)\varphi^2(s)\bigr)\varphi(a)\varphi^2(b)\quad
 \text{(as $s\tau(s)=1$)}\\
 &=s\varphi(a)\varphi^2(b)\quad\text{(as
 $s\varphi(s)\varphi^2(s)=1$)}\\
 &=s\tau\bigl((\tau\varphi(a))(\tau\varphi^2(b))\bigr)\\
 &=\overline{a\cdot b}.
\end{split}
\]
Now, because of \cite[Theorem 2.4]{EO}, for any $a,b,c\in A$,
\[
\begin{split}
\bar b\cdot (a\cdot c)
 &=(s\tau(b))\cdot\bigl((\tau\varphi(a))(\tau\varphi^2(b))\bigr)\\
 &=\bigl(\tau\varphi(s\tau(b))\bigr)
  \bigl(\tau\varphi^2(\tau(\varphi(a)\varphi^2(c)))\bigr)\\
  &=\bigl(\tau\varphi(s)\bigr)\varphi^2(b)
    \varphi\bigl(\varphi(a)\varphi^2(c)\bigr)\\
  &=\bigl(\tau\varphi(s)\bigr)\varphi^2(ab)c,
\end{split}
\]
which is symmetric on $a$ and $b$. Therefore,
\[
\delta_1(a,b)(c)=\bar b\cdot (a\cdot c)-\bar a\cdot(b\cdot c)=0
\]
and, similarly,
\[
\delta_2(a,b)(c)=(c\cdot a)\cdot\bar b -(c\cdot b)\cdot\bar a=0,
\]
for any $a,b,c\in A$. This shows that $(A,\cdot,\bar\ )$ is
trivially ($\delta=0$) a normal Lie related triple algebra.
\end{proof}

In our case, with $\calA=k[t,t^{-1},(1-t)^{-1}]$, the action of
$S_3$ is determined by the automorphisms $\varphi_\calA$ and
$\tau_\calA$ in Proposition \ref{th:automorphisms}, and the
distinguished element $s$ is $-t'$, which satisfies
$s\varphi_\calA(s)\varphi_\calA^2(s)=-t't''t=1=s\tau_\calA(s)$.

\bigskip

\section{The solution to \cite[Problem 13.3]{Paul}}

In this section, the ideals of both the Tetrahedron algebra and the
Onsager algebra will be determined. Let us first deal with the
Tetrahedron algebra, whose determination constitutes the problem
posed in \cite[Problem 13.3]{Paul}.

\begin{proposition}
$\frgtetra$ is a prime Lie algebra.
\end{proposition}
\begin{proof}
Since $\frsldos\otimes\calA\otimes_{\calA}k(t)$ is isomorphic to the
simple Lie algebra $\frsldos\otimes k(t)$, it follows that
$\frg=\frsldos\otimes\calA$ is prime, and hence so is $\frgtetra$.
\end{proof}

Let $\fra$ be an ideal of $\frg=\frsldos\otimes\calA$. The aim is to
prove that there is an ideal $\calI$ of $\calA$ such that
\[
\fra=u_0\calI\oplus u_1\calI\oplus u_2\calI.
\]
Therefore, the set of ideals of $\frg$ (and hence of $\frgtetra$) is
in bijection with the set of ideals of $\calA$.

To prove this, take any element $u_0a+u_1b+u_2c$ in $\frg$
($a,b,c\in\calA$). Then, using Theorem \ref{th:Abasis} and since
$tt't''=-1$:
\begin{equation}\label{eq:u0u1u2x}
\begin{split}
[u_0,[u_1,[u_2,u_0a+u_1b+u_2c]]]&=-u_1b,\\
[u_1,[u_2,[u_0,u_0a+u_1b+u_2c]]]&=-u_2c,\\
[u_2,[u_0,[u_1,u_0a+u_1b+u_2c]]]&=-u_0a.
\end{split}
\end{equation}

\smallskip

\begin{theorem} The ideals of $\frg=\frsldos\otimes\calA$ are
precisely the subspaces
\[
u_0\calI\oplus u_1\calI\oplus u_2\calI,
\]
where $\{u_0,u_1,u_2\}$ is the $\calA$-basis of $\frg$ in Theorem
\ref{th:Abasis} and $\calI$ is an ideal of $\calA$.
\end{theorem}
\begin{proof}
Let $\fra$ be an ideal of $\frg$ and consider the subspaces
$\calI_i=\{a\in\calA: u_ia\in\fra\}$, $i=0,1,2$. As in
\eqref{eq:u0u1u2x}, $[u_2b,[u_0,[u_1,u_0a]]]=-u_0ba$, so $\calI_0$
is an ideal of $\calA$, and so are $\calI_1$ and $\calI_2$. Now,
because of \eqref{eq:u0u1u2x}, $\fra=u_0\calI_0\oplus
u_1\calI_1\oplus u_2\calI_2$. But for any $a\in\calI_0$ and
$b\in\calA$, $[u_2b,u_0a]=-u_1bt''a\in\fra$, so $\calA t''a\subseteq
\calI_1$. Since $t''$ is invertible in $\calA$, this shows that
$\calI_0=\calA\calI_0=\calA t''\calI_0\subseteq \calI_1$. In the
same vein, one proves $\calI_0\subseteq \calI_1\subseteq
\calI_2\subseteq \calI_0$, so $\calI_0=\calI_1=\calI_2$ and the
result follows.
\end{proof}

Notice that $\calA$, as a ring of fractions of $k[t]$, is a
principal ideal domain.

\smallskip

Actually, a much more general result can be given:

\begin{theorem}\label{th:idealsgotimesA}
Let $\frg$ be a central simple finite dimensional Lie algebra, and
let $A$ be a unital commutative associative algebra. Then the ideals
of $\frg\otimes A$ are precisely the subspaces $\frg\otimes I$,
where $I$ is an ideal of $A$ ($\frg\otimes I$ is naturally
identified with a subspace of $\frg\otimes A$).
\end{theorem}
\begin{proof}
Let $\calM(\frg)$ be the associative subring (not necessarily
unital) of $\End_k(\frg)$ generated by $\{\ad_g:g\in\frg\}$ (see
\cite[Chapter X]{Jacobson}). A well-known result by Wedderburn shows
that $\calM(\frg)=\End_k(\frg)$, since $\frg$, as a central simple
Lie algebra, is an irreducible module for $\calM(\frg)$, and the
centralizer of the action is $k$. Let $\{g_i:i=1,\ldots,n\}$ be a
basis of $\frg$ and let $e_{ij}\in\End_k(\frg)$ be the linear map
given by $e_{ij}(g_k)=\delta_{jk}g_i$. Then, for any $x=\sum_{i=1}^n
g_i\otimes a_i\in\frg\otimes A$ ($a_i\in A$, $i=1,\ldots,n$),
$g_i\otimes a_i=(e_{ii}\otimes 1)(x)$ and $g_j\otimes
a_ib=(e_{ji}\otimes b)(x)$ belong to the ideal of $\frg\otimes A$
generated by $x$, as both $e_{ii}$ and $e_{ji}$ belong to
$\calM(\frg)$. Now, if $\fra$ is an ideal of $\frg\otimes A$, and
$I=\{ a\in A: g_1\otimes a\in \fra\}$, the above arguments show
that, for any $i=1,\ldots,n$, $I=\{a\in A: g_i\otimes a\in\fra\}$,
$I$ is an ideal of $A$, and $\fra=\frg\otimes I$.
\end{proof}

\smallskip

In \cite{DateRoan1}, an ideal $\fra$ of a Lie algebra $\frg$ is said
to be \emph{closed} if $Z(\fra)=\{x\in\frg: [x,\frg]\subseteq
\fra\}$ equals $\fra$. Note that $Z(\fra)/\fra$ is the center of
$\frg/\fra$, so the ideal $\fra$ is closed if the center of
$\frg/\fra$ is trivial.

\begin{corollary}
Let $\frg$ be a central simple finite dimensional Lie algebra, and
let $A$ be a unital commutative associative algebra. Then any ideal
of $\frg\otimes A$ is closed.
\end{corollary}
\begin{proof}
 From Theorem \ref{th:idealsgotimesA} it follows that any ideal of
$\frg\otimes A$ is of the form $\frg\otimes I$ for some ideal $I$ of
$A$, and hence the quotient $\frg\otimes A/\frg\otimes I$ is
isomorphic to $\frg\otimes (A/I)$, whose center is trivial, as $A/I$
is unital.
\end{proof}

\medskip

The closed ideals of the Onsager algebra have been determined in
\cite{DateRoan1} and \cite{DateRoan2}. Here all the ideals of the
Onsager algebra will be determined. To do so, let us identify the
Onsager algebra $\calO$ with $\psi(\Omega)=u_0(t-1)k[t]\oplus
u_1k[t]\oplus u_2tk[t]$ (Proposition \ref{pr:Omega}), which is
closed under the action of Klein's $4$ group. Consider the following
elements in $\calO$:
\[
v_0=u_0(t-1),\qquad v_1=u_1,\qquad v_2=u_2t,
\]
which are free generators of $\calO$ over $k[t]$, and satisfy
\begin{equation}\label{eq:vis}
[v_0,v_1]=-v_2(t-1),\qquad [v_1,v_2]=-v_0,\qquad [v_2,v_0]=v_1t,
\end{equation}
because of Theorem \ref{th:Abasis}, as $t't=t-1$ and $t''(t-1)=-1$.

Recall that the centroid of a Lie algebra $\frg$ over $k$ is the
centralizer of the adjoint action: $\Gamma(\frg)=\{
f\in\End_k(\frg): f([x,y])=[x,f(y)]\ \text{for any}\ x,y\in\frg\}$.

\begin{lemma}
\null\quad\null
\begin{romanenumerate}
\item $\calO$ is generated, as an algebra over $k$, by $v_0$, $v_1$
and $v_2$.
\item The centroid of $\calO$ is isomorphic to $k[t]$.
\item $\calO$ is prime.
\end{romanenumerate}
\end{lemma}
\begin{proof}
Since $(\ad_{v_2})^2(v_0t^n)=v_0t^{n+1}$, it follows that $v_0t^n$
belongs to the subalgebra generated by $v_0$, $v_1$ and $v_2$ for
any $n$. But then so does $v_1t^{n+1}=[v_2,v_0t^n]$ and
$v_2t^n(t-1)=[v_1,v_0t^n]$. Hence (i) follows.

It is clear that $k[t]$ embeds in the centroid $\Gamma(\calO)$,
since $\calO$ is an algebra over $k[t]$. But for any
$f\in\Gamma(\calO)$, $f(v_i)\in\{x\in\calO : [x,v_i]=0\}=v_ik[t]$,
so there are polynomials $p_i(t)\in k[t]$ such that
$f(v_i)=v_ip_i(t)$, $i=0,1,2$. Then
\[
v_0p_0(t)=f(v_0)=f([v_2,v_1])=\begin{cases}
 [f(v_2),v_1]=[v_2p_2(t),v_1]=v_0p_2(t),&\\
 [v_2,f(v_1)]=[v_2,v_1p_1(t)]=v_0p_1(t),&\end{cases}
\]
so $p_0(t)=p_1(t)=p_2(t)$. Because of (i), $f$ is determined by its
action on $v_0$, $v_1$ and $v_2$, and hence $f$ is the right
multiplication by $p_0(t)$.

Finally, $\calO\otimes_{k[t]} k(t)$ is a three dimensional simple
algebra over the field $k(t)$, so that $\calO$ is prime.
\end{proof}

\smallskip

In order to determine the ideals of $\calO$, first note that for any
ideal $J$ of $k[t]$, $\calO J$ is an ideal of $\calO$. These ideals
are closed under the action of Klein's $4$ group.

Let now $\calI$ be an ideal of $\calO$, and consider the following
subspace of $k[t]$:
\[
J_{\calI}=\{ p(t)\in k[t]: v_0p(t)+v_1p_1(t)+v_2p_2(t)\in \calI\
\text{for some}\ p_1(t),p_2(t)\in k[t]\}.
\]

\begin{proposition}\label{pr:idealsOidealskt}
Let $\calI$ be an ideal of $\calO$. Then:
\begin{romanenumerate}
\item $J_{\calI}$ is an ideal of $k[t]$.
\item $\calI$ lies between the ideals $\calO J_{\calI}t(t-1)$ and
$\calO J_{\calI}$:
\[
\calO J_{\calI}t(t-1)\subseteq \calI\subseteq
\calO J_{\calI}.
\]
\end{romanenumerate}
\end{proposition}
\begin{proof}
For any $p(t)\in J=J_{\calI}$, there are polynomials
$p_1(t),p_2(t)\in k[t]$ such that $x=v_0p(t)+v_1p_1(t)+v_2p_2(t)$
belongs to $\calI$. But
\[
(\ad_{v_2})^2(x)=v_0tp(t)+v_1tp_1(t)\in \calI,
\]
so $tp(t)\in J$, and hence $J$ is an ideal of $k[t]$.

Now, if $x=v_0p_0(t)+v_1p_1(t)+v_2p_2(t)\in\calI$, then $p_0(t)$
lies in $J$, but since
\[
\begin{split}
[v_2,x]&=v_0p_1(t)+v_1tp_0(t)\in\calI, \\
[-v_1,x]&=v_0p_2(t)-v_2p_0(t)(t-1)\in\calI,
\end{split}
\]
it follows that both $p_1(t)$ and $p_2(t)$ lie in $J$ too.
Therefore, $\calI$ is contained in $v_0J\oplus v_1J\oplus v_2J=\calO
J$. Moreover,
\[
\begin{split}
[v_0,[v_1,x]]&=v_1p_0(t)t(t-1)\in\calI,\\
[v_2,[v_0,[v_1,x]]]&=v_0p_0(t)t(t-1)\in\calI,\\
[-v_0,[v_2,x]]&=v_2p_0(t)t(t-1)\in\calI,
\end{split}
\]
so $\calO Jt(t-1)=v_0Jt(t-1)\oplus v_1Jt(t-1)\oplus v_2Jt(t-1)$ is
contained in $\calI$.
\end{proof}

\begin{corollary}
Let $\calI$ be a maximal ideal of $\calO$, then the quotient algebra
$\calO/\calI$ is either a one dimensional Lie algebra over $k$, or a
three dimensional simple Lie algebra over a finite field extension
of $k$.
\end{corollary}
\begin{proof}
By maximality, either $J_{\calI}=k[t]$, or $\calI=\calO J_{\calI}$
and $J=J_{\calI}$ is a maximal ideal of $k[t]$.  In the first case,
$v_0\in \calI$, and hence $v_2(t-1)=[v_1,v_0]\in\calI$ and
$v_1t=[v_2,v_0]\in \calI$. Thus $\calO/\calI$ is spanned by $\bar
v_1=v_1+\calI$ and $\bar v_2=v_2+\calI$, which satisfy $[\bar
v_1,\bar v_2]=0$. It follows that $\calO/\calI$ is abelian and, by
maximality of $\calI$, the dimension of $\calO/\calI$ is $1$.

Otherwise $\calI=\calO J_{\calI}$ and $J=J_{\calI}$ is a maximal
ideal of $k[t]$. Then there are three different possibilities,
according to $J$ being the ideal generated by $t$, by $t-1$, or by a
monic irreducible polynomial different from $t$ or $t-1$.

If $J=tk[t]$, then $\calO/\calI$ is spanned by $\bar v_i=v_i+\calI$,
$i=0,1,2$, which satisfy $[\bar v_2,\bar v_0]=0$, $[\bar v_1,\bar
v_2]=-\bar v_0$ and $[\bar v_0,\bar v_1]=\bar v_2$, thus giving a
three dimensional solvable Lie algebra, which contradicts the
maximality of $\calI$. The same happens if $J=(t-1)k[t]$.

However, if $J=p(t)k[t]$, for a monic irreducible polynomial
different from $t$ and $t-1$, then $K=k[t]/J$ is a finite field
extension of $k$, and $\calO/\calI$ is naturally a Lie algebra over
$K$ with a basis $\{\bar v_i= v_i+\calI: i=0,1,2\}$. Because of
\eqref{eq:vis}, this is a simple three dimensional Lie algebra over
$K$.
\end{proof}

Note that in $\calA=k[t,t^{-1},(1-t)^{-1}]$, both $t$ and $t-1$ are
invertible, and hence, with the same arguments as in the previous
proof, it is easily checked that the quotient of the Tetrahedron
algebra $\frgtetra$ by any maximal ideal is always a three
dimensional simple Lie algebra over a finite field extension of $k$.

\smallskip

Note that $k[t]$ decomposes as
\[
k[t]=t(t-1)k[t]\oplus \bigl( kt\oplus k(t-1)\bigr)
\]
(direct sum of subspaces). Therefore,
\begin{equation}\label{eq:OOtt1}
\calO=\calO t(t-1)\oplus \espan{v_it,v_i(t-1): i=0,1,2},
\end{equation}
and if $J=q(t)k[t]$ is a nonzero ideal of $k[t]$, then
\[
\calO J=\calO Jt(t-1)\oplus\espan{v_iq(t)t,v_iq(t)t(t-1): i=0,1,2}.
\]
Besides, for any nonzero ideal $J$ of $k[t]$, $[\calO t(t-1),\calO
J]$ is contained in $\calO Jt(t-1)$, so there are natural bijections
\[
\begin{matrix}
\bigl\{\text{ideals $\calI$ of $\calO$ with $\calO Jt(t-1)\subseteq
\calI\subseteq \calO J$}\bigr\}\\
\updownarrow\\
\bigl\{\text{$\calO$-submodules of  $\calO J/\calO Jt(t-1)$}\bigr\}\\
\updownarrow\\
\bigl\{\text{$\calO/\calO t(t-1)$-submodules of  $\calO J/\calO
Jt(t-1)$}\bigr\}
\end{matrix}
\]

Given an element $x\in\calO$ (respectively $x\in \calO J$), let us
denote by $\bar x$ its class modulo $\calO t(t-1)$ (respectively,
modulo $\calO Jt(t-1)$). Thus, from \eqref{eq:OOtt1}, with
$J=q(t)k[t]\ne 0$  and $w_i=v_iq(t)$, $i=0,1,2$,
\[
\begin{split}
\calO/\calO t(t-1)&=\bigoplus_{i=0}^2\bigl(k\overline{v_it}\oplus
k\overline{v_i(t-1)}\bigr),\\
\calO J/\calO Jt(t-1)&=\bigoplus_{i=0}^2\bigl(k\overline{
w_it}\oplus k\overline{w_i(t-1)}\bigr).
\end{split}
\]
Note also that for any $i,j$,
$[\overline{v_it},\overline{w_j(t-1)}]=0=[\overline{v_i(t-1)},\overline{w_jt}]$,
since $[v_it,w_j(t-1)],\,[v_i(t-1),w_jt]\in\calO Jt(t-1)$.

The eigenvalues of the action of $\overline{v_2t}$ on $\calO J/\calO
Jt(t-1)$ are:

\begin{itemize}
\item $0$, with eigenspace
$k\overline{w_2t}\oplus\bigl(\oplus_{i=0}^2k\overline{w_i(t-1)}\bigr)$,\\[1pt]

\item $1$, with eigenspace $k\overline{w_0t+w_1t}$, and\\[1pt]

\item $-1$, with eigenspace $k\overline{w_0t-w_1t}$.
\end{itemize}
(For instance, by \eqref{eq:vis}, $[v_2t,w_0t+w_1t]=w_1t^3+w_0t^3$,
but $t^3=t(t-1)(t+1)+t$, so $\overline{w_it^3}=\overline{w_it}$.)

Also, the eigenvalues of the action of $\overline{v_1(t-1)}$ on
$\calO J/\calO Jt(t-1)$ are:

\begin{itemize}
\item $0$, with eigenspace
$\bigl(\oplus_{i=0}^2k\overline{w_it}\bigr)\oplus
k\overline{w_1(t-1)}$,\\[1pt]

\item $1$, with eigenspace $k\overline{w_0(t-1)+w_2(t-1)}$,
and\\[1pt]

\item $-1$, with eigenspace $k\overline{w_0(t-1)-w_2(t-1)}$.
\end{itemize}

Besides, $[v_0t,w_1t]=-w_2t^2(t-1)$ belongs to $\calO Jt(t-1)$, and
the same happens to $[v_1t,w_0t]$, $[v_0(t-1),w_2(t-1)]$ and
$[v_2(t-1),w_0(t-1)]$. Moreover, $\overline{w_2t}$ generates the
$\calO$-submodule $\oplus_{i=0}^2k \overline{w_it}$, and
$\overline{w_1(t-1)}$ generates $\oplus_{i=0}^2
k\overline{w_i(t-1)}$.

Therefore, since any $\calO$-submodule of $\calO J/\calO Jt(t-1)$ is
the direct sum of its intersections with the previous eigenspaces,
we get:

\begin{proposition}\label{pr:idealsO}
Let $J=q(t)k[t]$ be a nonzero ideal of $k[t]$. Then the ideals
$\calI$ of $\calO$ with $J=J_{\calI}$ are the subspaces
\[
\calI=\calO Jt(t-1)\oplus \calS,
\]
where $\calS$ is of one of the following types:
\begin{romanenumerate}
\item $\calS=k\epsilon(w_0t+w_1t)\oplus k\delta(w_0t-w_1t)\oplus
k\epsilon\delta\gamma w_2t\oplus k\epsilon'(w_0(t-1)+w_2(t-1))\oplus
k\delta'(w_0(t-1)-w_2(t-1))\oplus k\epsilon'\delta'\gamma'
w_1(t-1)$, where $w_i=v_iq(t)$, $i=0,1,2$, and
$\epsilon,\delta,\gamma,\epsilon',\delta',\gamma'$ are either $0$ or
$1$, with $\epsilon+\delta\ne 0\ne\epsilon'+\delta'$ (as
$J=J_{\calI}$).
\item $\calS=\calS_\eta=\espan{
w_0t,w_1t,w_0(t-1),w_2(t-1),w_2t+\eta w_1(t-1)}$, with $0\ne \eta\in
k$.
\end{romanenumerate}
\end{proposition}

\begin{remark}
The ideals in Proposition \ref{pr:idealsO} which are invariant under
the action of Klein's $4$ group are those of type (i) with
$\epsilon=\delta=\epsilon'=\delta'=1$. In this case,
$\calI=v_0J\oplus v_1J_1\oplus v_2J_2$, where $J_1$ is either $J$ or
$Jt$, and $J_2$ is either $J$ or $J(t-1)$ (four possibilities which
correspond to $\gamma$ and $\gamma'$ being either $0$ or $1$).
\end{remark}

\smallskip

Recall that an ideal $\calI$ of $\calO$ is closed if $Z(\calI)=\{
x\in \calO: [x,\calO]\subseteq \calI\}$ equals $\calI$. Write
$J=J_{\calI}$. Then, for any $x=v_0p_0(t)+v_1p_1(t)+v_2p_2(t)$ in
$Z(\calI)$,
\[
\begin{split}
[x,v_1]&=-v_2(t-1)p_0(t)+v_0p_2(t)\in \calI\subseteq \calO J,\\
[x,v_2]&=-v_1tp_0(t)-v_0p_1(t)\in\calI\subseteq \calO J,
\end{split}
\]
and hence $p_1(t),p_2(t)\in J$ and $(t-1)p_0(t),tp_0(t)\in J$, so
$p_0(t)=tp_0(t)-(t-1)p_0(t)\in J$ too. Therefore,
\[
\calO Jt(t-1)\subseteq \calI\subseteq Z(\calI)\subseteq \calO J.
\]
Also, if $J=q(t)k[t]$ and $w_i=v_iq(t)$, $i=0,1,2$, using
\eqref{eq:OOtt1}, we get
\[
\begin{split}
[w_2t,\calO]&=[w_2t,\calO t(t-1)]+k[w_2t,v_0t]+k[w_2t,v_1t]\\
 &\hspace{2in} +\sum_{i=0}^2k[w_2t,v_i(t-1)]\\
 &\subseteq \calO Jt(t-1)+kw_1t+kw_0t,
\end{split}
\]
and this is contained in $\calI$ in case $\calI$ is as in item (ii)
of Proposition \ref{pr:idealsO}, so $w_2t\in Z(\calI)\setminus
\calI$ in this case. The same happens for $\calI$ as in item (i) of
Proposition \ref{pr:idealsO} with $\epsilon=\delta=1$ and
$\gamma=0$, or with $\epsilon'=\delta'=1$ and $\gamma'=0$. For the
remaining ideals in Proposition \ref{pr:idealsO}, it is easily
checked that they are closed. Therefore:

\begin{proposition}\label{pr:closedidealsO}
Let $\calI$ be a closed ideal of $\calO$, with $0\ne
J=J_{\calI}=q(t)k[t]$, then, with $w_i=v_ik[t]$, $i=0,1,2$, $\calI$
is one of  the following ideals:
\begin{alphaenumerate}
\item $\calI=\calO Jt(t-1)\oplus k\bigl(w_0t\pm w_1t\bigr)\oplus
k\bigl(w_0(t-1)\pm w_2(t-1)\bigr)$ ($4$ possibilities).
\item $\calI=\calO Jt(t-1)\oplus \bigl(\oplus_{i=0}^2kw_it\bigr)\oplus
k\bigl(w_0(t-1)\pm w_2(t-1)\bigr)$ ($2$ possibilities). In this case
$\calI =\calO Jt\oplus k(w_0\pm w_2)$.
\item $\calI=\calO Jt(t-1)\oplus k\bigl(w_0t\pm w_1t\bigr)\oplus
\bigl(\oplus_{i=0}^2 kw_i(t-1)\bigr)$ ($2$ possibilities). In this
case $\calI=\calO J(t-1)\oplus k(w_0\pm w_1)$.
\item $\calI=\calO J$.
\end{alphaenumerate}
\end{proposition}
\begin{proof}
Only the last assertions in (b) and (c) need to be checked. Since
$Jt=Jt(t-1)\oplus kq(t)t$, it follows that $\calO Jt=\calO
Jt(t-1)\oplus\bigl(\oplus_{i=0}^2kw_it \bigr)$. Besides,
$w_i(t-1)+w_i=w_it\in\calO Jt$. The last assertion in item (b)
follows at once. The argument for item (c) is similar.
\end{proof}

Note that only the ideal in (d) is invariant under the action of
Klein's $4$ group.

\bigskip

\section{Another presentation of the Tetrahedron algebra}

In this section, a new presentation of the Tetrahedron algebra,
based on the properties of our basis over $\calA$, will be given.

Consider the Lie algebra $\frf$ generated by the elements $z_0$,
$z_1$ and $z_2$ subject to the relations:
\begin{subequations}\label{eq:otherrel}
\begin{align}
&[[z_i,z_{i+1}],z_{i+2}]=0,\label{eq:otherrel1}\\
&[z_i,[z_i,z_{i+1}]]=z_{i+1}+[z_{i+2},z_i],\label{eq:otherrel2}\\
&\bigl[[z_{i+1},[z_{i+1},[z_{i+1},z_i]]],[z_{i+1},z_i]\bigr]=0,\label{eq:otherrel3}
\end{align}
\end{subequations}
for any $i\in\{0,1,2\}$, and where the indices are taken modulo $3$.

The aim of this section is to show that $\frf$ is isomorphic to the
Tetrahedron algebra $\frgtetra$, thus providing a different
presentation of this latter algebra.

First, let us obtain some consequences of the relations
\eqref{eq:otherrel}:

\begin{lemma}
For any $i,j\in\{0,1,2\}$, and taking indices modulo $3$, the
following relations are satisfied:
\begin{align}
&[(\ad_{z_i})^2(z_j),z_j]=0,\label{eq:otherrel4}\\
&\bigl[[z_{i-1},z_i],[z_i,z_{i+1}]\bigr]=z_i,\label{eq:otherrel5}\\
&(\ad_{[z_i,z_{i+1}]})^3(z_{i+1})
 =(\ad_{z_{i+1}})^2(z_i)-(\ad_{z_{i+1}})^4(z_i).\label{eq:otherrel6}
\end{align}
\end{lemma}
\begin{proof}
Equation \eqref{eq:otherrel4} is clear for $i=j$. Also,
\[
[(\ad_{z_i})^2(z_{i+1}),z_{i+1}]=[z_{i+1}+[z_{i+2},z_i],z_{i+1}]=0,
\]
because of \eqref{eq:otherrel2} and \eqref{eq:otherrel1}, while
\[
[(\ad_{z_i})^2(z_{i-1}),z_{i-1}]=[z_i,[[z_i,z_{i-1}],z_{i-1}]]=
 -[(\ad_{z_{i-1}})^2(z_i),z_i]=0
\]
by the previous argument. Hence \eqref{eq:otherrel4} holds.

Now,
\[
\begin{split}
[[z_{i-1},z_i],[z_i,z_{i+1}]]
 &=[z_{i-1},[z_i,[z_i,z_{i+1}]]]\qquad\text{(because of
 \eqref{eq:otherrel1})}\\
 &=[z_{i-1},z_{i+1}+[z_{i+2},z_i]]\qquad\text{(using
 \eqref{eq:otherrel2})}\\
 &=[z_{i-1},z_{i+1}]+[z_{i-1},[z_{i-1},z_i]]\\
 &=[z_{i-1},z_{i+1}]+z_i+[z_{i+1},z_{i-1}]
  \qquad\text{(using again
 \eqref{eq:otherrel2})}\\
 &=z_i.
\end{split}
\]
thus proving \eqref{eq:otherrel5}.

Finally,
\[
\begin{split}
(\ad_{[z_i,z_{i+1}]}&)^3(z_{i+1})\\
 &=\Bigl[[z_i,z_{i+1}],\bigl[[z_i,z_{i+1}],[[z_i,z_{i+1}],z_{i+1}]\bigr]\Bigr]\\
 &=-\Bigl[[z_i,z_{i+1}],\bigl[z_i,[[[z_i,z_{i+1}],z_{i+1}],z_{i+1}]\bigr]\Bigr]\quad\text{(by
 \eqref{eq:otherrel4})}\\
 &=\Bigl[[z_i,[z_i,z_{i+1}]],[[[z_i,z_{i+1}],z_{i+1}],z_{i+1}]\Bigr]\quad\text{(by
 \eqref{eq:otherrel3})}\\
 &=-\Bigl[z_{i+1}+[z_{i+2},z_i],[z_{i+1},[z_{i+1},[z_{i+1},z_i]]]\Bigr]\quad\text{(by
 \eqref{eq:otherrel2})}\\
 &=-(\ad_{z_{i+1}})^4(z_i)+\Bigl[z_{i+1},\bigl[z_{i+1},[[z_{i-1},z_i],[z_i,z_{i+1}]]\bigr]\Bigr]
  \quad\text{(by \eqref{eq:otherrel1})}\\
 &=-(\ad_{z_{i+1}})^4(z_i)+(\ad_{z_{i+1}})^2(z_i),\quad\text{(by
 \eqref{eq:otherrel5})}
\end{split}
\]
thus proving \eqref{eq:otherrel6}.
\end{proof}

\medskip

\begin{theorem}
There is a Lie algebra isomorphism $\Phi:\frgtetra\rightarrow \frf$
such that
\begin{equation}\label{eq:tetisof}
\begin{aligned}
\Phi(X_{01})&=2(z_2-[z_1,z_2]),&\qquad \Phi(X_{23})&=2(z_2+[z_1,z_2]),\\
\Phi(X_{02})&=2(z_0-[z_2,z_0]),& \Phi(X_{31})&=2(z_0+[z_2,z_0]),\\
\Phi(X_{03})&=2(z_1-[z_0,z_1]),& \Phi(X_{12})&=2(z_1+[z_0,z_1]).
\end{aligned}
\end{equation}
\end{theorem}
\begin{proof}
To check that the formulas above define a Lie algebra homomorphism,
it is enough to check the relations \eqref{eq:rel}. This is
straightforward. For instance, let us check that
$[\Phi(X_{01}),\Phi(X_{13})]=2\bigl(\Phi(X_{01})+\Phi(X_{13})\bigr)$
and that $[\Phi(X_{01}),[\Phi(X_{01}),[\Phi(X_{01}),\Phi(X_{23})]]]=
4[\Phi(X_{01}),\Phi(X_{23})]$. First,
\[
2\bigl(\Phi(X_{01})+\Phi(X_{13})\bigr)=
2\bigl(\Phi(X_{01})-\Phi(X_{31})\bigr)=-4\bigl([z_1,z_2]+[z_2,z_0]\bigr),
\]
while
\[
\begin{split}
[\Phi(X_{01})&,\Phi(X_{13})]\\
 &=-4\bigl[z_2-[z_1,z_2],z_0+[z_2,z_0]\bigr]\\
 &=-4\bigl([z_2,z_0]+[z_2,[z_2,z_0]]-\bigl[[z_1,z_2],[z_2,z_1]\bigr]\bigr)
  \quad\text{(using \eqref{eq:otherrel1})}\\
 &=-4\bigl([z_2,z_0]+z_2+[z_1,z_2]-z_2\bigr)\quad\text{(because of
 \eqref{eq:otherrel2} and \eqref{eq:otherrel5})}\\
 &=-4\bigl([z_2,z_0]+[z_1,z_2]\bigr).
\end{split}
\]
Now,
\[
\begin{split}
4[\Phi(X_{01}),\Phi(X_{23})]&=16\bigl[z_2-[z_1,z_2],z_0+[z_2,z_0]\bigr]\\
 &=-32[z_2,[z_2,z_1]]=-32(\ad_{z_2})^2(z_1),
\end{split}
\]
while
\[
\begin{split}
[\Phi(X_{01})&,[\Phi(X_{01}),[\Phi(X_{01}),\Phi(X_{23})]]]\\
 &=16\Bigl[z_2-[z_1,z_2],\bigl[z_2-[z_1,z_2],[z_2-[z_1,z_2],z_2+[z_1,z_2]]\bigr]\Bigr]\\
 &=-32\Bigl[z_2-[z_1,z_2],\bigl[z_2-[z_1,z_2],[z_2,[z_2,z_1]]\bigr]\Bigr]\\
 &=-32\Bigl((\ad_{z_2})^4(z_1)+(\ad_{[z_1,z_2]})^3(z_2\\
 &\qquad -\bigl[[z_1,z_2],[z_2,[z_2,[z_2,z_1]]]\bigr]-
    \bigl[z_2,[[z_1,z_2],[z_2,[z_2,z_1]]]\bigr]\Bigr)\\
 &=-32\Bigr((\ad_{z_2})^2(z_1)-2\bigl[[z_1,z_2],[z_2,[z_2,[z_2,z_1]]]\bigr]\Bigr)
  \quad\text{(by \eqref{eq:otherrel6})}\\
 &=-32(\ad_{z_2})^2(z_1)\quad\text{(because of \eqref{eq:otherrel3}).}
\end{split}
\]

Hence \eqref{eq:tetisof} defines a Lie algebra homomorphism $\Phi$,
which is onto, since the generators $z_0$, $z_1$ and $z_2$ of $\frf$
lie in the image of $\Phi$.

On the other hand, in $\frg=\frsldos\otimes\calA$, Theorem
\ref{th:Abasis} shows that,
\[
[u_0,[u_0,u_1]]=-[u_0,tu_2]=-tt''u_1=(1-t'')u_1=u_1+[u_2,u_0],
\]
and cyclically. Therefore, the relation \eqref{eq:otherrel2} is
satisfied by the $u_i$'s, and so are the relations
\eqref{eq:otherrel1} and \eqref{eq:otherrel3} because of the
$\bZ_2\times\bZ_2$-grading
$\frg=\frg_0\oplus\frg_1\oplus\frg_2=u_0\calA\oplus u_1\calA\oplus
u_2\calA$ (see \eqref{eq:tg0g1g2} and Theorem \ref{th:tg0g1g2}),
where $[\frg_i,\frg_{i+1}]=\frg_{i+2}$ and $[\frg_i,\frg_i]=0$ for
any $i$.

Hence, there is a surjective homomorphism $\Gamma: \frf\rightarrow
\frg$ such that $\Gamma(z_i)=u_i$, $i=0,1,2$. But the composition
$\Gamma\Phi$ satisfies
\[
\Gamma\Phi(X_{01}+X_{23})=2(u_2-[u_1,u_2])+
2(u_2+[u_1,u_2])=4u_2=\Psi(X_{01}+X_{23})
\]
(see \eqref{eq:Psi} and \eqref{eq:Abasis}), and also
$\Gamma\Phi(X_{02}+X_{31})=4u_0=\Psi(X_{02}+X_{31})$ and
$\Gamma\Phi(X_{03}+X_{12})=4u_1=\Psi(X_{03}+X_{12})$. Therefore,
$\Gamma\Phi$ coincides with $\Psi$ on a set of generators, and hence
$\Gamma\Phi$ coincides with the isomorphism $\Psi$. In particular,
$\Phi$ is one-to-one too.
\end{proof}

\medskip

\begin{corollary}
The Tetrahedron algebra $\frgtetra$ is isomorphic to the Lie algebra
generated by three elements $z_0$, $z_1$ and $z_2$, subject to the
relations \eqref{eq:otherrel}.
\end{corollary}

\bigskip

\section{$S_4$-action on the universal central extension of the
Tetrahedron algebra}

In \cite[Definition 3.3]{GeorgiaPaul}, a Lie algebra $\frgtetrahat$
is defined with generators
\[
\bigl\{\tilde X_{ij}:i,j\in\{0,1,2,3\},\, i\ne
j\bigr\}\cup\bigl\{C_p:p\in P\bigr\},
\]
where $P$ is the set of partitions of $\{0,1,2,3\}$ into two
disjoint subsets of two elements each, subject to the relations
\begin{romanenumerate}
\item $C_p$ is central for any $p\in P$,
\item $\sum_{p\in P} C_p=0$,
\item $\tilde X_{ij}+\tilde X_{ji}=C_{ij}$ for $i\ne j$, where
$C_{ij}=C_{p_{ij}}$, and $p_{ij}$ consists of $\{i,j\}$ and
$\{0,1,2,3\}\setminus\{i,j\}$,
\item $[\tilde X_{ij},\tilde X_{jk}]=2\bigl(\tilde X_{ij}+\tilde
X_{jk})$ for mutually distinct $i,j,k$ such that $(i,j,k)$ is even
(that is, the permutation $\sigma$ of $\{0,1,2,3\}$ with
$\sigma(0)=i$, $\sigma(1)=j$ and $\sigma(2)=k$ is even).
\item $[\tilde X_{hi},[\tilde X_{hi},[\tilde X_{hi},\tilde X_{jk}]]]
 =4[\tilde X_{hi},\tilde X_{jk}]$ for
mutually distinct $h,i,j,k$.
\end{romanenumerate}

This Lie algebra $\frgtetrahat$ is a central extension of
$\frgtetra$, and the kernel of the natural projection
$\pi:{\frgtetrahat}\rightarrow \frgtetra$ ($\tilde X_{ij}\mapsto
X_{ij}$, $C_p\mapsto 0$) is the two dimensional space spanned by
$\{C_p: p\in P\}$. Moreover, if the characteristic of the ground
field is $0$, then $\frgtetrahat$ is shown to be the universal
central extension of $\frgtetra$ \cite[Theorem 5.3]{GeorgiaPaul}.

The Lie algebra $\frgtetrahat$ presents a natural $A_4$-symmetry,
where $\sigma(\tilde X_{ij})=\tilde X_{\sigma(i)\sigma(j)}$, for any
$i\ne j$, and $\sigma(C_p)=C_{\sigma(p)}$ for any $p\in P$, where
$\sigma(p)$ is the partition obtained from $p$ by applying the
permutation $\sigma$ to its two components.

\smallskip

However, the automorphism group of any perfect Lie algebra embeds in
the automorphism group of its universal central extension (see
\cite[Proposition 1.3(v)]{vdK} or \cite[Proposition 2.2]{Arturo}).
Therefore, the Lie algebra $\frgtetrahat$ should show a symmetry
over the whole symmetric group $S_4$. Let us show how to modify
slightly the above generating set of $\frgtetrahat$ so as to make
clear this symmetry.

To do this, consider the new elements
\[
Y_{ij}=\tilde X_{ij}-\frac{1}{2}C_{ij}
\]
for distinct $i,j$. These elements satisfy
\[
Y_{ij}+Y_{ji}=\tilde X_{ij}+\tilde X_{ji}-C_{ij}=0
\]
by relation (iii) above, and for distinct $i,j,k$:
\[
\begin{split}
[Y_{ij},Y_{jk}]&=[\tilde X_{ij},\tilde X_{jk}]\\
&=\left\{\begin{aligned} 2\tilde X_{ij}+2\tilde X_{ji}&=
  2Y_{ij}+2Y_{ji}+C_{ij}+C_{jk}\\
   &=2Y_{ij}+2Y_{jk}-C_{ik}\quad\text{for even $(i,j,k)$}\\[2pt]
   2\tilde X_{ij}+2\tilde X_{ji}&+2C_{ik}=
    2Y_{ij}+2Y_{jk}+C_{ik}\quad\text{for odd
    $(i,j,k)$}\end{aligned}\right.
\end{split}
\]
(see \cite[Lemma 3.5]{GeorgiaPaul}). Hence,
\begin{equation}\label{eq:YYC}
[Y_{\sigma(0)\sigma(1)},Y_{\sigma(1)\sigma(2)}]=
 2Y_{\sigma(0)\sigma(1)}+2Y_{\sigma(1)\sigma(2)}
 -(-1)^\sigma C_{\sigma(0)\sigma(2)}
\end{equation}
for any $\sigma\in S_4$ ($(-1)^\sigma$ denotes the signature of
$\sigma$).

Therefore, the generating set $\bigl\{ Y_{ij}:i,j\in\{0,1,2,3\},\,
i\ne j\bigr\}\cup\bigl\{C_p:p\in P\bigr\}$ satisfies the relations
\begin{romanprimeenumerate}
\item $C_p$ is central for any $p\in P$,
\item $\sum_{p\in P} C_p=0$,
\item $Y_{ij}+Y_{ji}=0$ for $i\ne j$,
\item $[Y_{\sigma(0)\sigma(1)},Y_{\sigma(1)\sigma(2)}]=
 2Y_{\sigma(0)\sigma(1)}+2Y_{\sigma(1)\sigma(2)}
 -(-1)^\sigma C_{\sigma(0)\sigma(2)}$ for any $\sigma\in S_4$
\item $[Y_{hi},[Y_{hi},[Y_{hi},Y_{jk}]]]
 =4[Y_{hi},Y_{jk}]$ for
mutually distinct $h,i,j,k$.
\end{romanprimeenumerate}

Now, the whole $S_4$ acts on these generators by
$\sigma(Y_{ij})=Y_{\sigma(i)\sigma(j)}$, for $i\ne j$, and
$\sigma(C_p)=(-1)^\sigma C_{\sigma(p)}$ for any $p\in P$. The
relations (i)'--(v)' above are invariant under this action of $S_4$.
Therefore $S_4$ embeds in the automorphism group
$\Aut(\frgtetrahat)$.

Note that for $i\ne j$,
\begin{multline*}
\sigma(\tilde X_{ij})=\sigma\bigl(Y_{ij}+\frac{1}{2}C_{ij}\bigr)=
Y_{\sigma(i)\sigma(j)}+\frac{1}{2}C_{\sigma(i)\sigma(j)}\\
= \tilde
X_{\sigma(i)\sigma(j)}-\frac{1}{2}(1-(-1)^\sigma)C_{\sigma(i)\sigma(j)},
\end{multline*}
so that
\[
\sigma(\tilde X_{ij})=\begin{cases} \tilde
X_{\sigma(i)\sigma(j)}&\text{if
$\sigma$ is even,}\\
\tilde X_{\sigma(i)\sigma(j)}-C_{\sigma(i)\sigma(j)}&\text{if
$\sigma$ is odd.}\end{cases}
\]
Also, the kernel of the Lie algebra epimorphism
$\pi:\frgtetrahat\rightarrow \frgtetra$ is spanned by $\{ C_p:p\in
P\}$, and the $C_p$'s are fixed by the elements of Klein's $4$
group. Hence, as in \eqref{eq:tg0g1g2},
\[
\frgtetrahat=\frt\oplus(\frgtetrahat)_0\oplus(\frgtetrahat)_1\oplus(\frgtetrahat)_2,
\]
where $\frt=\espan{C_p:p\in P}$, and the restriction
$\pi\vert_{(\frgtetrahat)_i}:(\frgtetrahat)_i\rightarrow
(\frgtetra)_i$ is an isomorphism for any $i=0,1,2$, where
$(\frgtetra)_i=\Psi^{-1}(\frg_i)$. It follows easily from here that
the involution and the binary multiplication of the normal Lie
related triple algebra associated to this $S_4$-action on
$\frgtetrahat$ coincide (up to isomorphism) with the ones already
considered for $\frgtetra$ (actually for $\frg=\Psi(\frgtetra)$).
Also, the elements (compare to \eqref{eq:Abasis})
\[
\begin{split}
\hat
u_0&=\frac{1}{4}\bigl(Y_{02}+Y_{31}\bigr)=\frac{1}{4}\bigl(\tilde
X_{02}+\tilde X_{31}-C_{02}\bigr),\\
\hat
u_1&=\frac{1}{4}\bigl(Y_{03}+Y_{12}\bigr)=\frac{1}{4}\bigl(\tilde
X_{03}+\tilde X_{12}-C_{03}\bigr),\\
\hat
u_2&=\frac{1}{4}\bigl(Y_{01}+Y_{23}\bigr)=\frac{1}{4}\bigl(\tilde
X_{01}+\tilde X_{23}-C_{01}\bigr),
\end{split}
\]
project onto the generators $\Psi^{-1}(u_0)$, $\Psi^{-1}(u_1)$ and
$\Psi^{-1}(u_2)$ of $\frgtetra$ (Theorem \ref{th:Abasis} and
Corollary \ref{co:Main}), and hence, since $\frgtetrahat$ is
perfect, they are generators of $\frgtetrahat$ as a Lie algebra over
$k$.

\providecommand{\bysame}{\leavevmode\hbox
to3em{\hrulefill}\thinspace}
\providecommand{\MR}{\relax\ifhmode\unskip\space\fi MR }
\providecommand{\MRhref}[2]{%
  \href{http://www.ams.org/mathscinet-getitem?mr=#1}{#2}
} \providecommand{\href}[2]{#2}


\begin{thebibliography}{KMRT98}

\bibitem[BT]{GeorgiaPaul}
G.~Benkart and P.~Terwilliger, \emph{The universal central extension
of the three-point $\frsldos$ loop algebra}, preprint
\texttt{arXiv:RA/0512422}.

\bibitem[DR00]{DateRoan1}
E.~Date and S.S.~Roan, \emph{The algebraic structure of the Onsager
algebra}, Czech. J.~Phys. \textbf{50} (2000), no.~1, 34--44.

\bibitem[DR00']{DateRoan2}
\bysame, \emph{The structure of quotients of the Onsager algebra by
closed ideals}, J.~Phys. A: Math. Gen. \textbf{33} (2000),
3275--3296


\bibitem[EO05]{EO}
Alberto Elduque and Susumu Okubo, \emph{Lie algebras with
$S_4$-action and structurable algebras}, preprint
\texttt{arXiv:math.RA/0508558}.

\bibitem[HT05]{Paul}
B.~Hartwig and P.~Terwilliger, \emph{The Tetrahedron algebra, the
Onsager algebra, and the $\frsldos$ loop algebra}, preprint
\texttt{arXiv:math-ph/0511004}.

\bibitem[Jac79]{Jacobson}
N.~Jacobson, Lie algebras, Dover Publications, New York 1979.

\bibitem[Oku05]{Okubo}
S.~Okubo, \emph{Symmetric triality relations and structurable
algebras}, Linear Algebra Appl. \textbf{396} (2005), 189--222.

\bibitem[Ons44]{Onsager}
L.~Onsager, \emph{Crystal statistics. I. A two-dimensional model
with an order-disorder transition}, Phys. Rev. (2)
  \textbf{65} (1944), 117--149.

\bibitem[Pia02]{Arturo}
A.~Pianzola, \emph{Automorphisms of toroidal Lie algebras and their
central quotients}, J.~Algebra Appl. \textbf{1} (2002), no.~1,
113--121.

\bibitem[vdK73]{vdK}
W.L.J.~van der Kallen, \emph{Infinitesimally central extensions of
Chevalley groups}, Lecture Notes in Mathematics, Vol. 356,
Springer-Verlag, Berlin, 1973.

\end{thebibliography}
\end{document}